\documentclass[11pt]{article}
\usepackage{amsfonts}
\usepackage{amssymb}

\newtheorem{theorem}{Theorem}[section]
\newtheorem{proposition}[theorem]{Proposition}
\newtheorem{lemma}[theorem]{Lemma}
\newtheorem{corollary}[theorem]{Corollary}
\newtheorem{definition}{Definition}[section]
\newtheorem{remark}{Remark}[section]
\newtheorem{example}{Example}[section]
 \usepackage{amssymb}
\usepackage{epsfig}
\usepackage{amsmath}
\usepackage{amsfonts}

\newcommand{\R}{{{\Bbb R}}}

\newcommand{\N}{{{\Bbb N}}}

\newcommand{\Q}{{{\Bbb Q}}}

\def\qed{\hbox to 0pt{}\hfill$\rlap{$\sqcap$}\sqcup$\medbreak}

\title{Riemann integration via primitives for a new proof to the change of variable theorem}
\date{May, 2011}
\begin{document}
\maketitle

\begin{center}
Dedicated to Hyman Kestelman and Roy O. Davies on the 50th anniversary of the publication of their excellent papers on change of variable.
\end{center}

\vspace{1cm}

\begin{center}
{\large Rodrigo L\'opez Pouso \\
Departamento de
An\'alise Matem\'atica\\
Facultade de Matem\'aticas,\\Universidade de Santiago de Compostela, Campus Sur\\
15782 Santiago de
Compostela, Spain.
}
\end{center}

\begin{abstract}
We approach the Riemann integral via generalized primitives to give a new proof for a general result on change of variable originally proven by Kestelman and Davies. Our proof is similar to Kestelman's, but we hope readers will find it clearer thanks to the use of a new test for the Riemann integrability (which we introduce in this paper) along with some ingredients from some other more recent proofs  available in the literature. We also include a bibliographical review of related results and proofs. Although this paper emphasizes in the change of variable theorem, our contributions to the Riemann integration theory are of independent interest. For instance, we present a very simple proof to the fact that continuous functions are integrable which avoids the use of uniform continuity.
\end{abstract}

\section{Introduction}
Let $G:I=[a,b] \longrightarrow \R$ have a continuous derivative on the interval $I$, and let $f:G(I) \longrightarrow \R$ be continuous on $G(I)$.

If $G$ is not constant then we can use the chain rule and the Fundamental Theorem of Calculus to compute
\begin{equation}
\nonumber
 \left( \int_{G(a)}^{G(t)}{f(x) \, dx}\right)'=f(G(t)) \, G'(t) \quad (t \in [a,b]),
\end{equation}
and now it suffices to integrate between $a$ and $b$ to finish the nontrivial part of the proof of the following well--known theorem for the Riemann integral:

\begin{theorem}
\label{t1}
The change of variables formula
\begin{equation}
\label{ch}
\int_{G(a)}^{G(b)}{f(x) \, dx}=\int_a^b{f(G(t))G'(t) \, dt}
\end{equation}
is valid provided that $G:I=[a,b] \longrightarrow \R$ is continuously differentiable on $I$ and $f:G(I) \longrightarrow \R$ is continuous on $G(I)$.
\end{theorem}

Theorem \ref{t1} plays a fundamental role in elementary courses on Riemann integration: it is the key for evaluating exactly thousands of integrals. Moreover, Theorem \ref{t1} is stronger than it might appear at first glance. Notice that the substitution $G$ needs not be differentiable everywhere: it suffices that we can split the right--hand side in (\ref{ch}) as a finite sum of integrals each of which satisfies the conditions of Theorem \ref{t1}. In fact, as we will show, much more general results are known in the theory of Riemann integration.

\medbreak

Ch. J. de la Vall\'ee Poussin established a change of variable formula for Lebesgue integrals at the beginning of the 20th century, see \cite{val}. The availability of that result reduced later interests of the mathematical community in similar ones with Riemann integrals. Therefore, it is not a surprise that we have to wait until 1961 for the following general theorem, when it was proven by H. Kestelman \cite{kes} and R. O. Davies \cite{dav} in two consecutive papers in the same journal. From now on, and unless stated otherwise, integrability is to be understood in the Riemann sense.

\begin{theorem}[Kestelman's Theorem]
\label{tkes}
Assume that $g:I=[a,b] \longrightarrow \R$ is integrable on $I$ and let $G(t)=c+\int_a^t{g(s) \, ds}$ for all $t \in I$ and some $c \in \R$.

If $f:G(I) \longrightarrow \mathbb R$ is integrable on $G(I)$ then $(f \circ G) \, g$ is integrable on $I$ and
 \begin{equation}
\nonumber
\int_{G(a)}^{G(b)}{f(x) \, dx}=\int_a^b{f(G(t))g(t) \, dt}.
\end{equation}
\end{theorem}

Kestelman proved Theorem \ref{tkes} first, and then Davies found an elementary proof which, in particular, avoids the concept of null measure set, involved in Kestelman's original proof. Notice that, as pointed out by Kestelman, the assumptions in Theorem \ref{tkes} imply that the change of variable formula is valid with Lebesgue integrals (see, for instance, \cite{mcs, ser} or Theorem 6.95 on page 325 in \cite{str}). However, change of variable theorems for Lebesgue integrals give no information about the integrability of $(f \circ G) g$ in the Riemann sense, and therefore Kestelman's Theorem is neither a particular case to any of them nor a piece of completely romantic mathematics.

\medbreak

\medbreak

Kestelman's Theorem and Davies' elementary proof have been revised and rediscovered several times, see \cite{bag, cat, mar, nav, pop, pre, sar, ter, var} and our last section on bibliographical reviews. In this introduction we only point out the work done by Preiss and Uher in \cite{pre}, who established the following more general version of Theorem \ref{tkes}:

\begin{theorem}
\label{preuhe}
Let $g:I=[a,b] \longrightarrow \R$ be integrable on $I$, $c \in \R$, and $G(t)=c+\int_a^t{g(s) \, ds}$ for all $t \in I$.

 If $f:G(I)\longrightarrow \R$ is bounded on $G(I)$ then $f$ is integrable on $G(I)$ if and only if $(f\circ G)g$ is integrable on $I$ and, in that case,
 we have
 \begin{equation}
 \label{cvf}
 \int_{G(a)}^{G(b)}{f(x) \, dx}=\int_a^b{f(G(t))g(t) \, dt}.
 \end{equation}

\end{theorem}

 In this paper we present a new proof to Theorem \ref{preuhe} which keeps within the theory of Riemann integration and which we hope that readers may find more readable than others. Roughly, we adapt the classical proof described at the beginning of the introduction, in the sense we use the fact that the mapping
$$t \longmapsto \int_{G(a)}^{G(t)}{f(x) \, dx}$$
is a ``generalized" primitive of $(f \circ G) g$. 

\medbreak

\section{Riemann integrability through primitives}
Teachers always warn students not to make the following typical mistake: the definition of the Riemann integral of a function $f:[a,b] \longrightarrow \R$ is 
$$\int_a^b{f(x) \, dx}=F(b)-F(a), \quad \mbox{where $F$ is a primitive of $f$ ($F'=f$ in $[a,b]$).}$$
As we all know, integrable functions need not have primitives at all, so the previous formula is not adequate as a definition.

It is however true that we can define the Riemann integral by means of primitives, see \cite{tal}, although not primitives in the usual sense. In this section we propose an easy way to do it and we highlight some interesting consequences.

We are going to use the most basic concepts from the Riemann integration theory, which can be looked up in textbooks such as \cite{spi, tre}.

\bigbreak

To motivate our new definition of primitive, let $f:[a,b] \longrightarrow \R$ be bounded and consider the indefinite lower integral
\begin{equation}
\label{lint}
F_*(x)=\underline{\int_a^x}{f(y) \, dy} \quad (x \in [a,b]).
\end{equation}
For $x, y \in [a,b]$, $x<y$, we have (thanks to the additivity of the lower integral with respect to the intervals of integration) 
$$F_*(y)-F_*(x)=\underline{\int_x^y}{f(z) \, dz} ,$$
and, therefore,
$$\inf_{x \le z \le y} f(z) \le \dfrac{F_*(y)-F_*(x)}{y-x} \le \sup_{x \le z \le y} f(z).$$
Notice that the same property is fulfilled by the function
\begin{equation}
\label{uint}
F^*(x)=\overline{\int_a^x}{f(y) \, dy} \quad (x \in [a,b]).
\end{equation}
We have just shown that bounded functions have ``generalized primitives" according to the following definition (which does not involve derivatives).
\begin{definition}
\label{dgp}
Let $f:[a,b] \longrightarrow \R$ be a bounded function. 

A generalized primitive of $f$ is a function $F:[a,b] \longrightarrow \R$ such that for all $x,y \in [a,b]$, $x<y$, we have
\begin{equation}
\label{gp}
\inf_{x \le z \le y} f(z) \le \dfrac{F(y)-F(x)}{y-x} \le \sup_{x \le z \le y} f(z).
\end{equation}
\end{definition}

Generalized primitives need not be differentiable, but (\ref{gp}) readily yields some relations between their Dini derivatives
and the function $f$.

\begin{proposition}
\label{prodin}
Let $f:[a,b]\longrightarrow \R$ be bounded and let $F$ be a generalized primitive of $f$.

For every $x \in [a,b)$ we have
\begin{eqnarray}
\label{eqdi+}
\min\left\{f(x),\liminf_{y \to x^+}f(y)\right\} &\le& D_+f(x) \\
\nonumber&\le& D^+f(x)\le \max \left\{f(x),\limsup_{y\to x^+}f(y)\right\},
\end{eqnarray}
and for every $x \in (a,b]$ we have
\begin{eqnarray}
\label{eqdi-}
\min\left\{f(x),\liminf_{y \to x^-}f(y)\right\} &\le& D_-f(x)\\
\nonumber &\le& D^-f(x)\le \max \left\{f(x),\limsup_{y\to x^-}f(y)\right\}.
\end{eqnarray}
In particular, if $f$ is continuous at $x \in [a,b]$ then $F'(x)=f(x)$ (as usual, by $F'(a)$ we mean the right--hand derivative, and by $F'(b)$ the left--hand derivative).
\end{proposition}

\noindent
{\bf Proof.} We will only prove the left--hand inequality in (\ref{eqdi+}) because the remaining inequalities can be proven in analogous ways. Let $x \in [a,b)$ be fixed; we simply take limit inferior in the first inequality of (\ref{gp}) and we get the desired result:
\begin{align*}
D_+f(x)=\liminf_{y \to x^+}\dfrac{F(y)-F(x)}{y-x} & \ge \liminf_{y \to x^+}\inf_{x\le z \le y}f(z)= \lim_{y \to x^+}\inf_{x\le z \le y}f(z)\\
&= \lim_{y \to x^+}
\min\left\{f(x), \inf_{x < z \le y}f(z)\right\} \\
&=\min\left\{f(x),\liminf_{y \to x^+}f(y) \right\}.
\end{align*}
The existence of derivatives when $f$ is continuous is an immediate consequence of (\ref{eqdi+}) and (\ref{eqdi-}).
\qed

When $f$ is continuous generalized primitives are exactly the usual ones.

\begin{corollary}
\label{cogp}
If $f:[a,b]\longrightarrow \R$ is continuous in $[a,b]$ then every gene\-ra\-li\-zed primitive of $f$ is a primitive in usual sense and, conversely, primitives are generalized primitives.
\end{corollary}

\noindent
{\bf Proof.} Proposition \ref{prodin} ensures that generalized primitives satisfy $F'=f$ in $[a,b]$, so they are primitives in the usual sense.

Conversely, if $F$ is a primitive of $f$ then for every $x,y \in [a,b]$, $x<y$, the Mean Value Theorem guarantees the existence of some $z \in (x,y)$ such that
$$\dfrac{F(y)-F(x)}{y-x}=F'(z)=f(z),$$
which implies (\ref{gp}). \qed

Now we are in a position to establish our first test for the Riemann integrability in terms of generalized primitives. The Barrow's Rule for computing integrals using generalized primitives is included in the test.
\begin{theorem}
\label{car}
Lef $f:[a,b] \longrightarrow \R$ be bounded and $A \in \R$. The following two statements are equivalent:
\begin{enumerate}
\item The function $f$ is integrable on $[a,b]$ and $\int_a^b{f(x) \, dx}=A$;
\item Every generalized primitive of $f$ satisfies $F(b)-F(a)=A$.
\end{enumerate}
 \end{theorem}
\noindent
{\bf Proof.}  Let $F:[a,b] \longrightarrow \R$ be an arbitrary generalized primitive of $f$ and let $P=\{x_0,x_1,\dots,x_n\}$ be a partition of $[a,b]$.

By (\ref{gp}), we have
$$F(b)-F(a)=\sum_{k=1}^{n}{[F(x_k)-F(x_{k-1})]} \le U(f,P),$$
and, similarly, $F(b)-F(a) \ge L(f,P)$. Since $P$ was arbitrary and $f$ is integrable on $[a,b]$, we conclude that $F(b)-F(a)=\int_a^b{f(x) \, dx}=A$.

\bigbreak

Conversely, we take into account that $F_*$ and $F^*$, as defined in (\ref{lint}) and (\ref{uint}), are generalized primitives of $f$. Hence condition 2 ensures that
$$\underline{\int_a^b}{f(s) \, ds}=F_*(b)-F_*(a)=A=F^*(b)-F^*(a)=\overline{\int_a^b}{f(s) \, ds}.$$
\qed

Forgetting about computations, we have the following concise test for the Riemann integrability.

\begin{theorem}
\label{cocar}
A bounded function $f:[a,b] \longrightarrow \R$ is integrable if and only if any pair of its generalized primitives differ in a constant.
\end{theorem}

\noindent
{\bf Proof.} Assume that $f$ is integrable and let $F$ and $G$ be two of its generalized primitives. For each $x \in (a,b]$ we use Theorem \ref{car} in the interval $[a,x]$ to deduce that $F(x)-F(a)=G(x)-G(a)$ or, equivalently, that $F(x)-G(x)=F(a)-G(a)$.

Conversely, if $F$ and $G$ are two arbitrary primitives of $f$ the assumption ensures the existence of some constant $c \in \R$ such that $F(x)=c+G(x)$ for all $x \in [a,b]$. Hence
$F(b)-F(a)=G(b)-G(a),$ which implies that $f$ is integrable on $[a,b]$ by Theorem \ref{car}. \qed

When $f$ is integrable then its generalized primitives are exactly its indefinite integrals. More precisely we have the following corollary.

\begin{corollary}
\label{coin}
If $f:[a,b]\longrightarrow \R$ is integrable then $F:[a,b]\longrightarrow \R$ is a generalized pri\-mi\-tive of $f$ if and only if for each $x_0 \in [a,b]$ there exists some $c \in \R$ such that
$$F(x)=c+\int_{x_0}^x{f(y) \, dy} \quad (x \in [a,b]).$$
\end{corollary}

\noindent
{\bf Proof.} Notice that the indefinite integral
$$x \in [a,b]\longmapsto \int_{x_0}^x{f(y)\, dy}$$
is a generalized primitive of $f$ which, according to Theorem \ref{cocar}, is equal to any other generalized primitive up to adding some constant.\qed

\begin{remark}
In \cite[Theorem 1]{tho} Thomson solves completely the question of which functions $F$ are indefinite integrals of some unknown Riemann--integrable function $f$. Corollary \ref{coin} reveals that indefinite integrals are exactly generalized primitives, so we conclude, thanks to \cite[Theorem 1]{tho}, that $F$ is a generalized primitive of some unknown integrable function if and only if $F$ satisfies condition (3) in \cite{tho}, namely, if for all $\varepsilon >0$ a positive $\delta$ can be found so that
$$\sum_{i=1}^n\left| \dfrac{F(\xi_i)-F(x_{i-1})}{\xi_i-x_{i-1}}- \dfrac{F(x_i)-F(\xi_i')}{x_i-\xi_{i}'}\right|
(x_i-x_{i-1})<\varepsilon$$
for every subdivision $a=x_0<x_1<\cdots<x_n=b$ that is finer than $\delta$ and every choice of associated points $x_{i-1}<\xi_i \le \xi_i'<x_{i}$.
\end{remark}

We are going to illustrate the applicability of Theorem \ref{cocar} in the following two well--known situations.

\begin{example}
\label{excont} 
Continuous functions are integrable, and the usual way to prove it uses uniform continuity. Theorem \ref{cocar} yields an alternative and very easy proof. In particular, we avoid the uniform continuity (perhaps the following proof alone justifies the introduction of generalized primitives).

\medbreak

\noindent
{\bf Proof.} Let $f:[a,b]\longrightarrow \R$ be continuous on $[a,b]$ and let $F$ and $G$ be two of its generalized primitives. According to Corollary \ref{cogp}, $F$ and $G$ are primitives in the usual sense, so we have $(F-G)'=0$ in $[a,b]$ and then $F-G$ is constant. We conclude that $f$ is integrable on $[a,b]$ by virtue of Theorem \ref{cocar}. \qed
\end{example}

Our second example concerns a typical non--integrable bounded function.

\begin{example}
Let $f(x)=0$ for $x \in [0,1]\cap \Q$ and $f(x)=1$ for $x \in [0,1]\setminus \Q$.

The functions
$$F_{\lambda}(x)=\lambda x \quad (x \in [0,1]),$$
are generalized primitives of $f$ provided that $\lambda \in [0,1]$. Obviously, $F_1$ and $F_0$ do not differ in a constant, hence $f$ is not integrable on $[0,1]$.
\end{example}

Using Theorem \ref{cocar} to deduce integrability needs checking that generalized primitives differ in a constant. The following lemma is very useful for that task,
and it will be fundamental in our proof of Theorem \ref{preuhe} in the next section. For different proofs see \cite[Proposition 1]{kol}, \cite[Theorem 2]{sar} or \cite[Lemma 6.89]{str}.

\begin{lemma}
\label{lema1}
If $F$ is Lipschitz on $[a,b]$ and $F'(x)=0$ for almost all $x \in [a,b]$ then $F$ is constant on $[a,b]$.
\end{lemma}

The combination of Lemma \ref{lema1} and Theorem \ref{cocar} is powerful, and we show it by proving the sufficient part of Lebesgue's test for Riemann--integra\-bi\-lity, which ensures 
that integrable functions are exactly those bounded functions whose sets of discontinuity points have zero measure. Here is the precise statement and its short proof.
\begin{proposition}
\label{leb}
If $f:[a,b]\longrightarrow \R$ is bounded and continuous almost everywhere in $[a,b]$ then $f$ is integrable on $[a,b]$.
\end{proposition}

\noindent
{\bf Proof.} Let $F$ and $G$ be two generalized primitives of $f$. Then $F-G$ is a Lipschitz continuous function satisfying $(F-G)'=0$ almost everywhere. Lemma \ref{lema1} guarantees now that $F-G$ is constant.\qed

Finally, note that Lebesgue's test and Proposition \ref{prodin} guarantee that if $f$ is integrable and $F$ is one of its generalized primitives, then $F'(x)=f(x)$  for almost all $x \in [a,b]$. Conversely, we have the following result.

\begin{proposition}
\label{gpd}
If $f:[a,b]\longrightarrow \R$ is bounded, $F:[a,b]\longrightarrow \R$ is absolutely continuous, and $F'(x)=f(x)$ a.e. in $[a,b]$, then $F$ is a generalized primitive of $f$.
\end{proposition}

\noindent
{\bf Proof.} (This proof is the only part in this paper where integrals are to be understood in the Lebesgue sense. Notice that we will not use Proposition \ref{gpd} in the next section.) For $x,y
\in [a,b]$, $x<y$, the Fundamental Theorem of Calculus for the Lebesgue integral yields
$$F(y)-F(x)=\int_x^y{F'(z) \, dz}=\int_x^y{f(z) \, dz},$$
which implies (\ref{gp}). \qed

\begin{remark}
We have shown a way to introduce the Riemann integral from generalized primitives immediately after studying lower and upper Darboux integrals and we saw that it has some advantages.

We point out now another advantage: the definition of gene\-ra\-lized primitive can be adapted for other intervals than compact ones, and it could therefore be possible to unify the proper and improper theories of Riemann integration thanks to generalized primitives. We will not go further with this idea, because similar and deeper approaches have already been succesfully tried. For instance, Talvila \cite{tal} defines integrals by means of primitives in the weak or distributional sense (to which the adjective ``generalized" would fit better than to ours).
\end{remark}

\section{Change of variables}

This section is devoted to proving Theorem \ref{preuhe} following the ideas in \cite{kes, sar} combined with our theorems \ref{car} and \ref{cocar}. The proof due to Sarkhel and V\'yborn\'y in \cite{sar} is similar to Kestelman's and it consists of two steps, proving, loosely speaking, Theorem \ref{preuhe} for $G$ monotone first and using it to establish the general case. Here we work directly with the general case.

We alse use in our proof the well--known result which guarantees that Lispchitz continuous functions map null--measure sets into null--measure sets. See \cite[Lemma 1]{kes} or \cite[Lemma 6.87]{str}.

\bigbreak

\noindent
{\bf Proof of Theorem \ref{preuhe}.} Assume that $f$ is integrable on $G(I)$, let ${ H}={ H}(t)$ be an arbitrary generalized primitive of $(f \circ G)g$ and let
$${\cal F}(t)=\int_{G(a)}^{G(t)}{f(x) \, dx} \quad (t \in [a,b]).$$
If we prove that ${\cal F}-{ H}$ is a constant then we can conclude that every generalized primitive of $(f�\circ G)g$ is of the form ${ H}=c+{\cal F}$ for some constant $c \in \R$, and therefore, Theorem \ref{car} guarantees that $(f\circ G)g$ is integrable and  
$$\int_{a}^b{f(G(t))g(t) \, dt}={\cal F}(b)-{\cal F}(a)=\int_{G(a)}^{G(b)}{f(x) \, dx}.$$

To prove that ${\cal F}-{ H}$ is constant we note that ${\cal F}-{ H}$ is Lipschitz and we use Lemma \ref{lema1}, so it suffices to prove that $({\cal F}-{ H})'=0$ almost everywhere in $(a,b)$. To do it, we first notice that we can neglect all those points $t \in (a,b)$ where $g$ is discontinuous, because they form a null--measure set. Therefore we have only to study the subset of $(a,b)$ where $g$ is continuous. We split this set into 
$$A=\{ t \in (a,b) \,: \, \mbox{$g$ continuous at $t$ and $g(t)=0$}\}$$
and 
$$B=\{ t \in (a,b) \, : \, \mbox{$g$ continuous at $t$ and $g(t) \neq 0$}\}.$$

For all $t \in A$ the function $(f \circ G)g$ is continuous at $t$, so Proposition \ref{prodin} ensures that ${ H}'(t)=f(G(t))g(t)=0.$ In turn, for $t \in A$ and $h \neq 0$, $|h|$ sufficiently small, we have
$$\left| \dfrac{{\cal F}(t+h)-{\cal F}(t)}{h}\right|\le  \left| \dfrac{\int_t^{t+h}g(s) \, ds}{h}\right|\sup_{x \in G(I)}|f(x)|,$$
and then ${\cal F}'(t)=0$. Hence $({\cal F}-{ H})'=0$ everywhere in $A$.

We turn our attention now to the set $B$, which we first decompose as $B=C \cup D$, where
$$C=\{ t \in B \, : \, \mbox{$f$ is continuous at $G(t)$}\} \, \, \, \mbox{and}�\, \, \, D=B \setminus C.$$
For all $t \in C$ the function $(f�\circ G)g$ is continuous at $t$, so Proposition \ref{prodin} ensures ${ H}'(t)=f(G(t))g(t)$, and the chain rule gives ${\cal F}'(t)=f(G(t))g(t)$. We have proven that $({\cal F}-{ H})'=0$ everywhere in $C$.

Finally, we prove that $$D=\{t \in (a,b) \, : \, \mbox{$g$ continuous at $t$, $g(t)\neq 0$, $f$ discontinuous at $G(t)$}\}$$ is a null--measure set. To do it we mimic some ideas in Serrin and Varberg's proof of \cite[Theorem 1]{ser} and in Kestelman's proof of \cite[Lemma 2]{kes}.\footnote{We are more explicit on this in our review of Kestelman's paper in the next section.}

For each $t \in D$ we can find a sufficiently large $n \in \N$ such that $|g| >1/n$ and $g$ does not change sign in $[t-1/n,t+1/n] $, so $D$ is expressible as a countable union of sets of the form
$$D_{\mu}=\{t \in D \, : \, \mbox{$|g|>\mu$ and $g$ has constant sign in $[t-\mu,t+\mu] \cap I$}\} \, \, (\mu >0).$$
Now we prove that for each $\mu>0$ the set $D_{\mu}$ is a null--measure set, and it suffices to show that $D_{\mu} \cap J$ is a null--measure set for an arbitrary interval $J \subset I$ with length less than $\mu/2$. 

If $D_{\mu} \cap J$ is empty we are done, so assume that there is some $t \in D_{\mu} \cap J$. The definition of $D_{\mu}$ and the length of $J$ ensure that $|g| >\mu$ and $g$ does not change sign in $J$. Hence the restriction $G_{|J}$ is strictly monotone and it has an inverse $(G_{|J})^{-1}$ which is Lipschitz continuous (with Lipschitz constant $1/\mu$). Now we have $D_{\mu}\cap J=(G_{|J})^{-1}(G(D_{\mu}\cap J)),$ a null--measure set because $G(D_{\mu} \cap J)$ is contained in the set of discontinuity points of $f$ and $(G_{|J})^{-1}$ is Lipschitz continuous. The proof that $D$ is a null--measure set is finished. 

Hence $({\cal F}-{ H})'=0$ a.e. in $[a,b]$, and the first part of the proof is over.

\bigbreak

Conversely, assume that $(f \circ G)g$ is integrable and define 
$${\cal H}(t)=\int_a^t{f(G(s))g(s) \, ds} \quad (t \in [a,b]).$$
Now let $F$ be a generalized primitive of $f$ in the interval $G(I)$, and consider the composition ${\cal F}=F \circ G$.

Adjusting the previous arguments one can show that $({\cal F}-{\cal H})'=0$ a.e. in $[a,b]$ (Proving that $D$ is a null--measure set is easier in this case: it suffices to note that $D$ is contained in the set of discontinuity points of $(f \circ G)g$, which is a null--measure set\footnote{It seems that Kestelman and Davies dealt with the hardest part in Theorem \ref{preuhe}.}). Hence there is some $c \in \R$ such that 
\begin{equation}
\label{n2}
{\cal F}(t)=F(G(t))=c+{\cal H}(t)\, \, \, \mbox{for all $t \in I$,}
\end{equation}
and, in particular,\begin{equation}
\label{n1}
\int_a^b{f(G(s))g(s) \, ds}={\cal H}(b)-{\cal H}(a)={\cal F}(b)-{\cal F}(a).
\end{equation}
Since $F$ was an arbitrary generalized primitive of $f$, we deduce from (\ref{n2}) that all primitives of $f$ are equal up to an additive constant, and therefore $f$ is integrable in the interval $G(I)$. Moreover we have, by virtue of Corollary \ref{coin}, that ${\cal F}(t)=F(G(t))=\hat c+\int_{G(a)}^{G(t)}{f(x) \, dx}$ for some $\hat c \in \R$, and then (\ref{n1}) implies (\ref{cvf}). \qed

\section{Bibliographical reviews on change of variable}
\noindent
(1959) Marcus establishes in \cite{mar} sufficient conditions on $G$ which imply that the composition $f \circ G$ is integrable. Subsequently, the author gives some results on change of variable for the Riemann integral. Marcus' paper leans on some previous works by Zaslavski on Riemann integrability of compositions which go back to 1953. We drawn all the previous information from Mathematical Reviews.

A remarkable example in \cite{kes} shows that, in the conditions of Kestelman's Theorem, the composition $f \circ G$ needs not be Riemann integrable. Therefore, Kestelman's Theorem generalizes the results in \cite{mar}.

\medbreak

\noindent
(1961) Kestelman and Davies publish their respective proofs of Theorem \ref{tkes} in \cite{dav, kes}. One can deduce from \cite{dav, kes} that both authors were aware of each other's work before the publication of their papers. 

We have followed Kestelman's proof in this paper, so the reader already has an idea about the arguments in it. Interested readers might find it useful to know that Lemma 2 in \cite{kes} is proven there with some inacuracies (the set $S_+$ is expressed as an union of null--measure sets, but not necessarily a countable union, so one cannot deduce that $S_+$ is null). Lemma 2 in \cite{kes} is however correct, as it is a particular case to Theorem 1 in \cite{ser}. Anyway, we have benefitted from the proof of \cite[Lemma 2]{kes}, because the ideas in it about Lipschitzian inverses helped us to simplify the arguments we needed from the proof of \cite[Theorem 1]{ser} and, in particular, we have not used Lebesgue's outer measure.

Modestly, Davies claims in \cite{dav} that his paper is ``essentially based on the same ideas", which is not true at all. It is based on the following test for the Riemann integrability: $f(x)$ is integrable over $[a,b]$ if and only if it is bounded and given any $\varepsilon, \, \eta >0$ there exists a subdivision of  $[a,b]$ such that the intervals in which the oscillation of $f$ is greater than $\eta$ have total length less than $\varepsilon$. Davies' elementary proof of Theorem \ref{tkes} is really a masterpiece: it simply uses the most basic elements in Riemann integration (namely, Riemann sums) in a wonderfully clever way. It is however very subtle and, therefore, hard to follow or to convey to students.

\medbreak

\noindent
(1968) Varberg \cite{var} extends Marcus' results by stating and proving the following theorem:

\begin{theorem}
\label{vart}
Let $f$ be Riemann integrable on $[a,b]$ and let $G$ be absolutely continuous on $[c,d]$ with $G([c,d])\subset[a,b]$ and $G'$ (up to redefinition on a null--measure set) Riemann integrable on $[c,d]$ . Then $(f\circ G)\cdot G'$ is Riemann integrable on $[c,d]$ and $\int_{G(c)}^{G(d)}f(x)\,dx=\int_c^df(G(t))G'(t)\,dt$.
\end{theorem}

Note that the assumptions on $G$ imply that $G(t)=G(a)+\int_a^t{G'(s)ds}$ for all $t \in [c,d]$, with $G'$ Riemann integrable. In turn, functions $G$ in the conditions of Kestelman's Theorem are absolutely continuous (Lipschitz continuous, actually) and $G'=g$ almost everywhere, so Varberg's Theorem is equivalent to Kestelman's.

\begin{remark}
Cleverly, Varberg imposes assumptions over $f$ on an interval $[a,b]$ containing the range of $G$, and not exactly on the range of $G$ as we did in our statement of Kestelman's Theorem.

Varberg's statement avoids the inconsistency of Riemann integrability of $f$ on a singleton set, which does arise in our statement in case $G$ is a constant function.

Why didn't we follow Varberg's more rigorous writing?

First, the change of variables formula really needs no assumption outside the range of $G$, and this should remain clear from the very beginning.

Second, in the well--known convention $\int_a^a{f(x) \, dx}=0$ ($a \in \R$) for the Riemann integral underlies the idea that ``all functions are Riemann integrable on singletons". We are convinced that most readers will understand from our statement that no hypothesis on $f$ is required when $G$ is constant, and this is exactly what we mean.

Finally, the case $G$ constant is trivial and no one would ever use a theorem to study it. In our opinion, statements should not be overloaded just to avoid a minor abuse of language that only occurs in a trivial situation.
\end{remark}

\medbreak

\noindent
(1970) Preiss and Uher proved Theorem \ref{preuhe} in \cite{pre} by adjusting Davies' ideas in \cite{dav}.
\medbreak

\noindent
(1981) Navr\'atil revises Preiss and Uher proof of Theorem \ref{preuhe} in \cite{nav}.

\medbreak

\noindent
(1985) Ter\"ehkina \cite{ter} proves the following result on change of variables:
\begin{theorem}
\label{tere}
Let $f$ be Riemann integrable on $[a,b]$ and let $G$ have a Riemann integrable derivative on $[c,d]$, $G([c,d])\subset[a,b]$, and $G(c)=a$ and $G(d)=b$. Then $(f\circ G)\cdot G'$ is Riemann integrable on $[c,d]$ and $\int_{G(c)}^{G(d)}f(x)\,dx=\int_c^df(G(t))G'(t)\,dt$.
\end{theorem}

Theorem \ref{tere} is a particular case to Kestelman's Theorem. Indeed, the Second Fundamental Theorem of Calculus (see \cite[page 286]{spi}) guarantees that $G(t)=G(c)+\int_c^t{G'(s) \, ds}$ for all $t \in [c,d]$.

Ter\"ekhina also introduces in \cite{ter} an example showing that $f \circ G$ needs not be Riemann integrable in the assumptions of Theorem \ref{tere}.

\medbreak

\noindent
(1996) Sarkhel and V\'yborn\'y \cite{sar} rediscover Theorem \ref{preuhe}, without being aware of \cite{dav, kes, nav, pre}. Their proof differs from those in \cite{dav, nav, pre} in the sense that they use Lemma \ref{lema1}. In this sense, the proof in \cite{sar} is closer to that in \cite{kes}. The proof given in section 3 follows essentially the steps in \cite{sar}, but we  avoid the use of the change of variable formula when $G$ is monotone.

Besides a new proof, we owe to \cite{sar} at least two very interesting remarks: first, Theorem \ref{preuhe} is better than its analogues for the Lebesgue integral in the sense that we can deduce Riemann integrability of $f$ from that of $(f\circ G) g$, and this is not possible in general with Lebesgue integrals\footnote{We know, however, that such a result for the Lebesgue integral is valid when $G$ is monotone. We are grateful to Professor Rudolf V\'yborn\'y for having sent us a manuscript of his with a proof.}; second, a nice example which shows that we cannot simply omit the assumption that $f$ be bounded.

\medbreak
\noindent
(1997) Popovici and Bencze rediscover Kestelman's Theorem in \cite{pop} and introduce a new elementary proof, similar to Davies'.\footnote{We are grateful to Professor Mih\'aly Bencze for having sent us the paper \cite{pop}.}

\medbreak
\noindent
(1998) Cater \cite{cat} studies the change of variable formula in terms of the Dini derivatives of $G$. The use of Dini derivatives in connection with the formula of change of variable is classical, see \cite{hob}, but Cater's approach is new and not restricted to monotone substitutions.

For completeness, we state the main result in \cite{cat}. Notice that our next statement gathers the information given in Theorem 1 and Corollary 2 in \cite{cat}:

\begin{theorem}
\label{tcat1}
Let $DG$ denote one of the four Dini derivatives of a continuous real valued function $G:I=[a,b]\longrightarrow \R$ (the same Dini derivative at all $x$). Let $G(a) \le G(b)$ and let $DG$ be bounded on $I$. Let $f$ be a bounded function on $G(I)$ such that for almost every $t \in [a,b]$, one or both of the functions $f \circ G$ or $DG$ is continuous at $t$.

In the previous conditions the upper and lower Darboux integrals satisfy the following inequality:
\begin{align*}
\overline{\int_a^b}f(G(t))\, DG(t) \, dt & \ge \overline{\int_{G(a)}^{G(b)}}f(x) \, dx \\
& \ge \underline{\int_{G(a)}^{G(b)}}f(x) \, dx \ge \underline{\int_a^b}f(G(t))\, DG(t) \, dt.
\end{align*}
Moreover, if $(f\circ G)\, DG$ is Riemann integrable on $I$, then $f$ is Riemann integrable on $G(I)$ and $\int_{G(a)}^{G(b)}f(x)\,dx=\int_a^b f(G(t))DG(t)\,dt$.
\end{theorem}

We know, see \cite[page 455]{hob}, that if $G$ satisfies the conditions of Theorem \ref{preuhe} and $G(a) \le G(b)$ then $G$ satisfies the conditions of Theorem \ref{tcat1}. In this sense, Cater's Theorem is better than Theorem \ref{preuhe} as it allows more types of substitutions. The price for that generality turns out to be the assumption on almost everywhere continuity of $f \circ G$ or $DG$, which cannot be omitted, see \cite[Example 3]{cat}.

On the other hand, we note that if $DG$ in Theorem \ref{tcat1} is Riemann integrable then, see \cite[page 456]{hob}, $G(t)=G(a)+\int_a^t{DG(s) \, ds}$ for all $t \in I$, so $G$ satisfies the assumptions in Theorem \ref{preuhe} with $g=DG$. Examples of functions $G$ in the conditions of Theorem \ref{tcat1} and such that $DG$ is not Riemann integrable are subtle, see \cite[Example 2]{cat}.

\medbreak
\noindent
(2001) Bagby wrote in \cite{bag} the most recent paper we known about Kestelman's Theorem. Bagby's version applies for functions $f$ assuming values in an arbitrary Banach space. The proof is elementary, similar to that in \cite{dav}.

\medbreak
\noindent
(2008) An extension of Kestelman's Theorem to higher dimensions with injective $G$ is introduced in \cite{mol}.

\end{document}